\documentclass[12pt,reqno]{article}
\usepackage{graphicx}
\usepackage[all]{xy}
\usepackage{amsmath}
\usepackage{amsthm}
\usepackage{latexsym,bm}
\usepackage{amssymb}
\usepackage{indentfirst}
\usepackage{float}

\usepackage{graphicx}%
\usepackage{multirow}%
\usepackage{amsmath,amssymb,amsfonts}%
\usepackage{amsthm}%
\usepackage{mathrsfs}%
\usepackage[title]{appendix}%
\usepackage{xcolor}%
\usepackage{textcomp}%
\usepackage{manyfoot}%
\usepackage{booktabs}%
\usepackage{algorithm}%
\usepackage{algorithmicx}%
\usepackage{algpseudocode}%
\usepackage{listings}%
\usepackage{enumerate}
\usepackage{setspace}
\usepackage{lineno}
\newtheorem{thm}{Theorem}[section]
\newtheorem{cor}[thm]{Corollary}
\newtheorem{lem}[thm]{Lemma}
\newtheorem{prop}[thm]{Proposition}
\newtheorem{defn}[thm]{Definition}
\newtheorem{rem}[thm]{Remark}

\numberwithin{equation}{section}

\begin{document}
	
	\title{\bf{\Large Bounds of Scalar curvature, S-curvature and distortion on $\infty$-Einstein Finsler manifolds}}
	
	\author{{\bf Bin Shen}}
	
	\date{}
	
	\maketitle
	
	\begin{quote}
		\small {\bf Abstract}. This manuscript investigates the curvature and topological properties of certain $\infty$-Einstein Finsler metrics on Finsler metric measure spaces. By imposing symmetry conditions, we construct a series of special metrics and analyze their equivalence on special manifolds. Provided a Ricci curvature bound, we establish a linear growth lower bound estimate for the S-curvature and the distortion, revealing the interplay between curvature and measure on $\infty$-Einstein Finsler manifolds.
		Furthermore, by introducing scalar curvature and imposing a linear growth lower bound condition, we derive upper and lower bounds for the distortion, S-curvature, and the scalar curvature itself on asymmetric essential gradient Ricci solitons with certain non-Riemannian curvature constraints. These results yield direct topological finiteness conclusions for some forward-complete $\infty$-Einstein Finsler manifolds. Our work partially addresses Gromov’s conjecture of scalar curvature in the context of Finsler metric measure spaces and provides a foundation for further research in geometric analysis within general Finsler geometry.
	\end{quote}

	\begin{quote}
	\end{quote}
	
	\baselineskip 17pt
	

	
\section{Introduction}
Ricci curvature is the most fundamental geometric concept in Finsler geometry. In \cite{Chern1996}, S.-S. Chern posed the following pivotal question:\\

\textit{Can every smooth manifold admit a Finsler metric with constant Ricci scalar, or at least one whose Ricci scalar is independent of the direction $y$?} \\

This problem is equivalent to finding an Einstein Finsler metric with constant or scalar factor on a given manifold. By imposing special metric structures or additional symmetry conditions, several classification theorems and characterization results have been established \cite{CST2012, LCC2020, Vill2023}. Notably, in \cite{LCC2020, Vill2023} and other works, the Einstein Finsler metrics obtained share the property of having  \emph{constant} Einstein factors.

Over the past decade, numerous studies have focused on Finsler metrics and curvatures in pursuit of the Chern conjecture. However, the influence of curvature properties on manifold structure remains incompletely understood. Unlike Riemannian geometry, a Finsler metric does not uniquely determine a canonical volume form on the manifold. Instead, a measure must be explicitly specified to form a Finsler metric measure space $(M, F, \mu)$. Common choices for volume forms include the Busemann-Hausdorff volume, the Holmes-Thompson volume, or more generally, any Borel measure. To quantify the deviation of the chosen measure from a ``canonical" one, two non-Riemannian curvatures--the distortion and the S-curvature--are introduced \cite{Shen1997}. In this framework, weighted Riemannian spaces emerge as a special subclass of Finsler metric measure spaces.

In 2009, S. Ohta defined weighted Ricci curvatures (see Sect. 3) on $(M,F,\mu)$, generalizing the Riemannian weighted Ricci curvature and corresponding to the curvature-dimension condition ($CD(K, N)$ condition) on general metric measure spaces \cite{Ohta2009}. These curvatures play a crucial role in geometric and analytic problems on manifolds \cite{Xia2014, Yin2019, Ohta2021}.
In this manuscript, we focus on the weighted Ricci curvature $Ric^{\infty}$, defined as
$$Ric^{\infty} (x,y):=Ric(x,y)+\dot{S}(x,y),$$
where $\dot{S}(x,y)$ denotes the derivative of the S-curvature along the geodesic emanating from $x$ in direction $y$. A Finsler metric measure space is called an \emph{$\infty$-Einstein Finsler manifold} if it admits that
\begin{eqnarray}\label{eq-inftyEF}
	Ric^{\infty}(x,y)=\sigma(x,y) F^2(x,y),
\end{eqnarray}
where $\sigma$ is a function on the sphere bundle. This concept is a special case of the broader \emph{weak $(a,b)$ weighted Einstein metric} introduced in \cite{ShenZhao2022}, defined by 
$$Ric_{a,b} = (n-1) \left(\frac{3\theta}{F}+\sigma\right)F^2,$$ 
where $Ric_{a,b}=Ric+a\dot S-b S^2$ generalizes both the weighted Ricci curvature \cite{Ohta2009,Ohta2021} and the projectively invariant weighted Ricci curvature \cite{ShenSun2021}.

Drawing inspiration from Perelman’s framework \cite{Pere2003}, we establish the following theorem, which describes how the Ricci curvature of an $\infty$-Einstein Finsler metric influences the growth of the S-curvature and distortion.
\begin{thm}\label{thm-Ricbound}
	Let $(M,F,p)$ be a forward complete $\infty$-Einstein Finsler manifold with a pole $p$, the Einstein factor $\sigma(x,y)\geq\frac12$ and the bounded Ricci curvature $|Ric|\leq cF^2$, for some $c>0$. Then for any point $x$ on $M$, it satisfies that
	\begin{align*}
		&S(x,y)\geq\frac12(d(p,x)-K_0),\\
		&\tau(x,y)\geq \frac{1}{4}(d(p,x)-K_0)^2-K'_0,
	\end{align*}
	where $K_0$ and $K_0'$ are two constants only depending on the dimension $n$, the Ricci curvature bound $c$, the Finsler metric $F$ and the measures on $S_pM$ and $S_xM$.\\  
\end{thm}

A more challenging problem in Riemannian geometry is understanding the geometric and topological implications of scalar curvature. Gromov proposed several related questions \cite{Grom2019}, including\\

\textit{What are the topologies of spaces of metrics, and the geometries of individual manifolds, with scalar curvature bounded from below?}\\

We extend these questions to the setting of metric measure spaces, particularly examining whether scalar curvature retains its influence in the Finslerian setting. However, compared to the extensive results on Ricci curvature, research on scalar curvature in Finsler manifolds remains rather limited.
In 1988, H. Akbar-Zadeh introduced a definition of scalar curvature in Finsler geometry via
$$r=\frac12g^{ij}Ric_{y^iy^j}, \quad\mbox{with}\quad Ric=R^k_{\,\,k}, $$
and established a generalized Schur theorem \cite{AZ1988}. However, this definition incorporates substantial Finslerian structure, making it too rigid for applications in mathematical and physical theories.
In this work, we propose an alternative scalar curvature defined by
\begin{eqnarray}\label{eq-R}
	R=g^{ij}R^{\,\,k}_{i\,\,kj},
\end{eqnarray}
which enhances geometric applicability by incorporating additional symmetry.
When the Einstein factor $\sigma(x)$ in \eqref{eq-inftyEF} is constant, the metric is called a \emph{gradient Ricci soliton}, The geometric properties of such Finsler metrics have recently been studied by Li-Mo-Wang \cite{LMW2024}, Q. Xia \cite{Xia2024}, and others. Gradient Ricci solitons clarify the interplay between curvature and measure.
In this manuscript, we introduce several specialized $\infty$-Einstein Finsler metrics by imposing further symmetry conditions. A particularly suitable candidate for deeper study is the \emph{asymmetric essential $\infty$-gradient Ricci soliton}, defined by 
\begin{eqnarray}
	Ric^{\infty}_y(V,W)=\sigma g_y(V,W),
\end{eqnarray}
for all vector fields $V,W$ on $M$, where $\sigma$ is a constant, which we may assume is $\sigma=\frac12$ after rescaling. 
Using lower bounds on the proposed scalar curvature \eqref{eq-R}, we derive estimates for the S-curvature, distortion, and scalar curvature itself on asymmetric essential $\infty$-gradient Ricci solitons. It shows the possibility of replacing the Ricci curvature by the scalar curvature on some special Finsler metric measure spaces to control the behaviour of curvature and measure.
\begin{thm}\label{thm-1.2}
	Let $(M,F,p)$ be a forward complete asymmetric essential Finsler gradient Ricci soliton with a pole $p$, the factor $\sigma=\frac12$ and 
	whose scalar curvature admits at least linear growth, i.e., $R\geq \gamma d(p,x)-\alpha$ for some positive constants $\alpha$ and $\gamma$. Suppose further that the curvature bound $|(C^{t}_{il}L^{li}_{s}-C^{t\,\,|i}_{is})R_{\,\,t}^{s}|\leq \frac{(n+1)}{2}K_1F$ holds for some $K_1\leq (n+1)\gamma$. Then, the distortion, S-curvature, and scalar curvature admit the following bounds
	\begin{align*}
		\frac14\left[d(p,x)-C_2\right]^2-\alpha-\beta\leq \tau(x,y)&\leq\frac14\left[d(p,x)+C_1\right]^2-\alpha-\beta-\gamma,\\
		|S|&\leq\left[\frac12 d(p,x)+C_1\right]F,\\
		\frac14\left[d(p,x)+C_3\right]^2-\alpha\leq R&\leq \frac14\left[d(p,x)+C_1\right]^2-\alpha,
	\end{align*}
	where $C_1$, $C_2$ and $C_3$ are constants only depending on $n$, $\gamma$ and the Finsler metric $F$ and the volume form on the local tangent sphere bundle $SB_p(1)$.
\end{thm}	
The non-Riemannian curvature condition $|(C^{t}_{il}L^{li}_{s}-C^{t\,\,|i}_{is})R_{\,\,t}^{s}|\leq \frac{(n+1)}{2}K_1F$ is mild, as it automatically holds in any precompact region (and hence any bounded domain) due to homogeneity. Such conditions are common in curvature comparison formulas (e.g., \cite{ShenShen2016}). Moreover, any compact Finsler gradient Ricci soliton (examples of which can be found in \cite{Xia2024}) naturally satisfies this condition, making it a standard assumption in non-Riemannian geometry. Moreover, Our discussion of growth of the scalar curvature $R=R(x,y)$ refers to behavior along each geodesic ray (with $y$ fixed in direction), ensuring the validity of the linear growth lower bound. \\

We anticipate that the concepts of the new scalar curvature and the $\infty$-Einstein metric, including its special cases, will play pivotal roles in subsequent studies of Finsler geometry, particularly in global geometric properties related to the analysis and topology. As a direct application of Morse theory, Theorems \ref{thm-Ricbound} and \ref{thm-1.2} yield the following topological classification for forward complete Finsler gradient Ricci solitons.
\begin{cor}
	A forward complete $\infty$-Einstein Finsler manifold with $\sigma(x,y)\geq\frac12$ admits finite topological types provided either the following
	\begin{itemize}
		\item  its Ricci curvature has a constant bounds, i.e., $|Ric|\leq cF^2$;
		\item  it is an asymmetric essential Finsler gradient Ricci soliton ($\sigma=\frac12$) with a scalar curvature satisfying the linear growth lower bound $R\geq \gamma d(p,x)-\alpha$ for positive constants $\alpha$ and $\gamma$, and the non-Riemannian curvature is bounded by $|(C^{t}_{il}L^{li}_{s}-C^{t\,\,|i}_{is})R_{\,\,t}^{s}|\leq \frac{(n+1)}{2}K_1F$, where $K_1\leq (n+1)\gamma$.
	\end{itemize}
\end{cor}
These results elucidate the relationship between curvature and measure on $\infty$-Einstein Finsler manifolds and provide a foundation for further topological investigations. This work not partially extends Perelman’s celebrated results but also aligns with Gromov’s profound forecasts. A natural follow-up question is whether the linear growth condition on scalar curvature can be further relaxed.\\

This manuscript is organized as follows. Section 2 provides a brief review of foundational concepts in Finsler geometry. In Section 3, we introduce new curvature tensors and $\infty$-Einstein Finsler metrics, along with their interrelations. The proof of Theorem \ref{thm-Ricbound} is presented in Section 4. Section 5 derives a key formula relating scalar curvature and distortion on a specific class of Finsler gradient Ricci solitons. Section 6 explores direct consequences for some special cases of Finsler metrics. Finally, in Section 7, we establish bounds for the scalar curvature, distortion, and S-curvature under lower bound assumption of scalar curvature  on asymmetric Finsler gradient Ricci soliton.

\section{Basic concepts on Finsler manifolds}
In this section, we review the fundamental concepts of Finsler geometry that will be used throughout this article. We always denote $SM$ to be the sphere bundle, $TM$ the tangent bundle, $TM_0:=TM\setminus\{0\}$ the punched tangent bundle, $T_xM$ the tangent space at $x$, etc.

A \emph{Finsler metric} $F$ on a differential manifold $M$ is a function $F:TM\rightarrow [0,+\infty)$ that assigns a norm to each tangent space. Specifically, $F$ satisfies
\begin{itemize}
	\item[(i)] \textbf{Smoothness}: $F$ is smooth and positive on $TM_0$;
	\item[(ii)] \textbf{Positive homogeneity}: $F(x,ky)=kF(x,y)$ for any $(x,y)\in TM$ and $k>0$;
	\item[(iii)] \textbf{Strong convexity}: For any $(x,y)\in TM_{0}$, the fundamental tensor matrix
	\begin{eqnarray}\label{Fiii}
		g_{ij}(x,y):=\frac12\frac{\partial F^2}{\partial y^i\partial y^j}(x,y)
	\end{eqnarray}
	is positive definite.
\end{itemize}
This generalizes the Riemannian metric, as the fundamental tensor $g$ depends on the tangent direction $y$, making the Finsler metric non-quadratic in $y$. The deviation from Riemannian geometry is represented by the \emph{Cartan tensor}.
\begin{eqnarray*}
	C(X,Y,Z):=C_{ijk}X^iY^jZ^k=\frac{1}{4}\frac{\partial^3F^2(x,y)}{\partial y^i\partial y^j\partial y^k}X^iY^jZ^k.
\end{eqnarray*}
This symmetric tensor has a trace called the \emph{mean Cartan tensor} $I=I_idx^i$, where $I_i=g^{jk}C_{ijk}$. 

Geodesics on a Finsler manifold satisfy the differential equation
\begin{eqnarray}\label{eq-geod}
	\frac{d^2x^i}{dt^2}+2G^i(x,\frac{dx}{dt})=0,
\end{eqnarray}
where the \emph{spray coefficients} $G^i$ are derived from the Finsler metric by 
$$G^i=\frac{1}{4}g^{ij}\left\{[F^2]_{y^jx^k}y^k-[F^2]_{x^j}\right\}.$$
Equation \eqref{eq-geod} is a second-order quasilinear ODE with unique local solutions that depend continuously on initial conditions $(x,y)$. Unlike Riemannian geometry, Finsler geodesics are generally \emph{irreversible}: a geodesic from $x_1$ to $x_2$ may not coincide with the reverse path. We therefore distinguish that
\begin{itemize}
	\item \emph{Forward geodesic}: a geodesic with initial condition $(x_1,y_1)\in S_{x_1}M$;
	\item \emph{Backword geodesic}: a geodesic with initial condition $(x_1,-y_1)\in S_{x_1}M$.
\end{itemize}
A Finsler metric is \emph{reversible} if $F(x,-y)=F(x,y)$ for any $(x,y)\in TM_0$, otherwise, it is called \emph{nonreversible}. The reversibilities
$$\rho_x=\sup_{y\in S_xM}\frac{F(x,-y)}{F(x,y)},\quad\mbox{and}\quad \rho_M=\sup_{x\in M}\rho_x,$$
quantify the asymmetry between forward and backward distances. On compact manifolds, $\rho_M$ is always finite.

The \emph{nonlinear connection coefficients} $N^i_j=\frac{\partial G^i}{\partial y^j}$ induce a horizontal-vertical decomposition of $TTM_0$, with horizontal basis $\frac{\delta}{\delta x^i}=\frac{\partial}{\partial x^i}-N^j_i\frac{\partial}{\partial y^j}$. The \emph{Chern connection} $\nabla$ is the unique torsion-free, almost metric-compatible connection on the puback bundle, satisfying
\begin{align*}
	\nabla_{u}v-\nabla_{v}u&=[u,v],\\
	w(\langle u,v\rangle_{y})-\langle\nabla_{w}u,v\rangle_{y}-\langle u,\nabla_{w}v\rangle_{y}&=2C_{y}(\nabla_{w}y,u,v),
\end{align*}
for any $u,v,w\in T_0M$, where $C$ is denote as the Cartan tensor. Its coefficients are locally given by 
$$\Gamma^i_{jk}=\frac12g^{il}\left\{\frac{\delta g_{jl}}{\delta x^k}+\frac{\delta g_{lk}}{\delta x^j}-\frac{\delta g_{jk}}{\delta x^l}\right\}.$$
A Finsler metric is \emph{Berwald} if $\Gamma^i_{jk}$ is $y$-independent (equivalently, if $G^i$ is quadratic in $y$). The Chern connection extends to $TM$ with covariant derivatives
\begin{itemize}
	\item \emph{Horizontal}: $T^i_{j|k}
	=\frac{\delta T^i_j}{\delta x^k}+\Gamma^i_{ll}T^l_j-\Gamma^l_{kj}T^i_l,$;
	\item \emph{Vertical}: $T^i_{j;k}=\frac{\partial T^i_j}{\partial y^k}$,
\end{itemize}
for some 2-tensor $T=T^i_j\frac{\delta}{\delta x^i}\otimes dx^j$ as an example.

For any smooth vector fields $X,W,Z$, the Chern curvature decomposes as
\begin{eqnarray}\label{Cherncurform}
	\Omega(X,W)Z=R(X,W)Z+P(X,\nabla_W(y),Z),
\end{eqnarray}
where the \emph{Chern-Riemann ciurvature} $R$ has components $R_j{}^i{}_{kl}=\frac{\delta\Gamma_{jl}^{i}}{\delta x^{k}}-\frac{\delta\Gamma_{jk}^{i}}{\delta x^{l}}+\Gamma_{km}^{i}\Gamma_{jl}^{m}-\Gamma_{lm}^{i}\Gamma_{jk}^{m},$ and the \emph{Chern non-Riemannian curvature} $P$ has $P_j{}^i{}_{kl}=-\frac{\partial\Gamma_{jk}^i}{\partial y^l}$.


The \emph{Landsberg curvature} $L:=L^i_{jk}\frac{\partial}{\partial x^i}\otimes dx^j\otimes dx^k$, where
\begin{eqnarray}
	L^i_{jk}:=[G^i]_{y^jy^k}-\Gamma_{jk}^i=-y^lP^{\,\,i}_{l\,\,jk}=C^i_{jk|0}=C^i_{jk|l}y^l,
\end{eqnarray}
measures the rate of change of the Cartan tensor along geodesics. Its trace
$$J=J_i dx^i, \quad J_i=g^{jk}L_{ijk}.$$
is the \emph{mean Landsberg curvature}. 
A Finsler metric is Berwald if and only if Chern non-Riemannian curvature vanishes, is Landsberg if $L=0$ and is  \emph{weak Landsberg metric} if $J=0$.

For further relations between these curvatures, we refer to \cite{BCS2000} and \cite{ShenShen2016}. Additional identities used in this work will be introduced in subsequent proofs.

\section{Curvatures and the $\infty$-Einstein Finsler metrics}
In this section, we introduce several well-known curvature notions along with new concepts explored in this manuscript, particularly a sequence of $\infty$-Einstein Finsler metrics and the scalar curvature.

The Chern Riemannian curvature gives rise to several important curvature quantities. First is the flag curvature, which generalizes the sectional curvature in Riemannian geometry. For any point $(x,y)\in TM_0$, and a 2-dimensional section $\Pi_y=\mbox{span}\{y,v\}$ in $T_xM$, we define the \emph{flag curvature} with pole $y$ as 
\begin{eqnarray}
	K(\Pi_y):=\frac{-R_{ijkl}(x,y)y^iv^jy^kv^l}{(g_{ik}(x,y)g_{jl}(x,y)-g_{il}(x,y)g_{jk}(x,y))y^iv^jy^kv^l}.
\end{eqnarray} 
Setting
\begin{eqnarray}
	R^i{}_k:=y^jR_{j\,kl}^{\,\,i}y^l,\quad R_{jk}:=g_{ij}R^i{}_k=y^iR_{ijkl}y^l=-R_{ijlk}y^iy^l=-R_{jikl}y^iy^l,
\end{eqnarray}
we obtain
$$K(\Pi_y)(v):=F^{-2}R_{jk}v^jv^k,$$
provided $g_{ij}(y)y^iv^j=0$. The quatities $R^i{}_k$ or $R_{jk}$ are called the flag curvature tensor. The \emph{Ricci curvature} is defined as
\begin{eqnarray}
	Ric:=y^iR_{i\,\,kl}^{\,k}y^l=g^{jk}R_{jk}=R^k_{\,\,k},
\end{eqnarray}
which is also referred to as the Ricci scalar since it represents a scalar function on $TM_0$.

As Finsler manifolds lack a canonical volume form, we define a Finsler metric measure space $(M,F,\mu)$ as a Finsler manifold $(M,F)$ equipped with a specified measure $\mu$, where $d\mu=\sigma(x)dx$. A notable example is the \emph{Busemann-Hausdorff volume form}
\begin{eqnarray}
	dV_{F}:=\sigma_{B}(x)\omega^{1}\wedge\cdots\wedge\omega^{n},
\end{eqnarray}
where
$$\sigma_{B}(x):=\frac{\mathrm{Vol}(B^{n}(1))}{\mathrm{Vol}\left(y=y^{i}\frac{\partial}{\partial x^i}\in T_{x}M|F(x,y)<1\right)},$$
with $\mathrm{Vol}$ denoting the Euclidean volume and $B^n(1)$ the Euledian unit ball.

The \emph{distortion} of a Finsler metric measure space is defined by 
$$\tau=\log\frac{\det(g_{ij})}{\sigma(x)},$$
which gives rise to the \emph{S-curvature} as 
$$S(x,y):=\frac{d}{dt}\left[\tau\left(\gamma(t),\dot{\gamma}(t)\right)\right]|_{t=0},$$
where $\gamma(t)$ is a geodesic with $\gamma(0)=x$ and $\dot\gamma (0)=y$.

Combining with the Ricci curvature and the S-curvature,  \cite{Ohta2021} define the weighted Ricci curvature on a Finsler metric measure space as follows for $N \in [n,+\infty]$.
\begin{itemize}
	\item $Ric^n(x,y):= \begin{cases}Ric(x,y)+\dot{S}(x,y) & \text { if }{S}(x,y)=0; \\ -\infty & \text { if }{S}(x,y)\neq 0,\end{cases}$
	\item $Ric^N(x,y):=Ric(x,y)+\dot{S}(x,y) - \frac{S^2(x,y)}{N-n}$ when $n<N<\infty$,
	\item $Ric^{\infty} (x,y):=Ric(x,y)+\dot{S}(x,y)$,
\end{itemize}
where $\dot{S}(x,y)$ denotes the derivative along the geodesic from $x$ in direction $y$. 

An $\infty$-Einstein Finsler metric is equivalently a $(1,0)$-weighted Einstein metric. 
Due to Schur's theorem in Riemannian geometry, Einstein metrics can choose the constant Einstein factor when the dimension of the manifold satisfies $n\geq 3$. However, this theorem does not always hold in Finsler geometry, which leads to some confusion about the name of the Einstein-Finsler metric/manifold. For example,
in prior works of Einstein-Finsler metrics \cite{ShenSun2021,ShenZhao2022,Xia2024}, $\sigma$ is assumed to be a function on $M$, while in \cite{GGPT2024,Vill2023,LCC2020}, $\sigma$ is taken as a 
constant. According to Theorem 1.1 in \cite{Xia2024}, the $\infty$-Einstein Finsler manifold is just the gradient almost Ricci soliton on a compact Finsler metric measure space. So we give the definitions as the follows.
\begin{defn}
	A Finsler metric measure space $(M,F,\mu)$ is called an $\infty$-Einstein Finsler manifold if it satisfies
	\begin{eqnarray}\label{eq-EF}
		Ric^{\infty}(x,y)=\sigma F^2(x,y),
	\end{eqnarray}
	where $\sigma(x,y)$ is a scalar function defined on $SM$. Moreover, $(M,F,\mu)$ is called a gradient almost Ricci solition if $\sigma(x)$ is a function on $M$, and is called a gradient Ricci soliton if $\sigma$ is a constant.
\end{defn}

While this definition is natural, its generality limits the scope for global geometric and analytic techniques—beyond Ricci curvature lower bounds as in Theorem \ref{thm-Ricbound}. To enable finer analysis, we introduce the following refined notions.
\begin{defn}
	A Finsler metric measure space $(M,F,\mu)$ is called an asymmetric $\infty$-Einstein Finsler manifold (resp. asymmetric gradient almost Ricci soliton or asymmetric gradient Ricci soliton) if
	\begin{eqnarray}
		Ric^{\infty}_y(y,V)=\sigma g_y(y,V),
	\end{eqnarray}
	for all vector fields $V$ on $M$, where $\sigma=\sigma(x,y)$ is a function on $SM$ (resp. $\sigma=\sigma(x)$ is a function on $M$ or $\sigma$ is a constant).
\end{defn}

\begin{rem}
	Expressed in coordinates as $Ric^{\infty}_y(y,V)=y^i(R^{\,\,k}_{j\,\,kl}+\tau_{|j|l})V^l$, this tensor is asymmetric: $Ric^{\infty}_y(y,V)\neq Ric^{\infty}_y(V,y)$. , justifying the term ``asymmetric".
\end{rem}

\begin{defn}
	A Finsler metric measure space $(M,F,\mu)$ is called a symmetric $\infty$-Einstein Finsler manifold (resp. symmetric gradient almost Ricci soliton or symmetric gradient Ricci soliton) if
	\begin{eqnarray}
		\frac12[Ric^{\infty}_y(y,V)+Ric^{\infty}_y(V,y)]=\sigma g_y(y,V),
	\end{eqnarray}
	for all vector fields $V$ on $M$, where $\sigma=\sigma(x,y)$ is a function on $SM$ (resp. $\sigma=\sigma(x)$ is a function on $M$ or $\sigma$ is a constant).
\end{defn}

We further propose an alternative symmetric Einstein-type condition.
\begin{defn}\label{def-VV}
	A Finsler metric measure space $(M,F,\mu)$ is called an essential $\infty$-Einstein Finsler manifold (resp. essential gradient almost Ricci soliton or essential gradient Ricci soliton) if
	\begin{eqnarray}
		Ric^{\infty}_y(V,V)=\sigma g_y(V,V),
	\end{eqnarray}
	for all vector fields $V$ on $M$, where $\sigma=\sigma(x,y)$ is a function on $SM$ (resp. $\sigma=\sigma(x)$ is a function on $M$ or $\sigma$ is a constant).
\end{defn}
Its asymmetric counterpart is 
\begin{defn}\label{def-asymFE}
	A Finsler metric measure space $(M,F,\mu)$ is called an asymmetric essential $\infty$-Einstein Finsler manifold (resp. asymmetric essential gradient almost Ricci soliton or asymmetric essential gradient Ricci soliton) if
	\begin{eqnarray}
		Ric^{\infty}_y(V,W)=\sigma g_y(V,W),
	\end{eqnarray}
	for all vector fields $V,W$ on $M$, where $\sigma=\sigma(x,y)$ is a function on $SM$ (resp. $\sigma=\sigma(x)$ is a function on $M$ or $\sigma$ is a constant).
\end{defn}
Analogous definitions apply to $(a,b)$-weighted Einstein Finsler manifolds. 

\begin{rem}
	When $\sigma$ is constant, metric rescaling allows us to normalize $\sigma$ to any preferred value without altering its sign. Henceforth, we fix $\sigma=\frac12$ for all Finsler gradient Ricci soliton in this manuscript.
\end{rem}

All known Einstein metrics in Finsler geometry contain abundant information about non-Riemannian curvatures, which describe the connections between geometry, physics, and analysis. However, this richness of non-Riemannian data can also obscure the key relationships on the spherical bundle. Definitions \ref{def-VV} and \ref{def-asymFE} strip away these extraneous non-Riemannian features, which is precisely why we refer to them as ``essential."
Moreover, in certain special cases, the above Einstein Finsler metrics are equivalent—further illustrating the significance of the term ``essential."

\begin{prop} 
	For a Landsberg Finsler metric, the following identities hold.
	\begin{itemize}
		\item The asymmetric essential gradient almost Ricci soliton coincides with the asymmetric gradient almost Ricci soliton.
		\item The gradient almost Ricci soliton is equivalent to the symmetric gradient almost Ricci soliton.
	\end{itemize}
\end{prop}

\begin{proof}
	We demonstrate the equivalences using the asymmetric $\infty$-Einstein Finsler metric as an example. In local natural coordinates, this metric satisfies 
	\begin{eqnarray}\label{new1.1}
		(R^{\,\,k}_{i\,\,kj}+\tau_{|i|j})y^i=\sigma g_{ij}y^i,
	\end{eqnarray}	
	in local natural coordinate system. 
	From the formula 
	\begin{eqnarray}
		R^{\,\,k}_{j\,\,kl;m}=P^{\,\,k}_{j\,\,km|l}-P^{\,\,k}_{j\,\,lm|k}-P^{\,\,k}_{j\,\,ks}L^s_{lm}+P^{\,\,k}_{j\,\,ls}L^s_{km},
	\end{eqnarray}
	we can find that
	\begin{eqnarray}\label{a-1}
		y^jR^{\,\,k}_{j\,\,km;h}=-L^k_{kh|m}+L^k_{mh|k}+L^k_{ks}L^s_{mh}-L^k_{ms}L^s_{kh}.
	\end{eqnarray}
	On the other hand, considering $[\frac{\delta}{\delta x^i},\frac{\partial }{\partial y^j}]=(\Gamma^k_{ij}+L^k_{ij})\frac{\partial}{\partial y^k}$, we derive 
	\begin{eqnarray}
		\tau_{|i|j;m}=\tau_{|i;m|j}-L^k_{jm}\tau_{|i;k}+P^{\,\,k}_{i\,\,jm}\tau_{|k},
	\end{eqnarray}
	and thus
	\begin{eqnarray}\label{a-2}
		\tau_{|i|j;m}y^i=(\tau_{|i;m}y^i)_{|j}-L^k_{jm}\tau_{|i;k}y^i-L^{k}_{jm}\tau_{|k}.
	\end{eqnarray}	
	Taking the vertical Chern covariant derivative of \eqref{new1.1}, we obtain 
	\begin{eqnarray}
		R^{\,\,k}_{i\,\,kj}+\tau_{|i|j}=\sigma g_{ij},
	\end{eqnarray}	
	on a Landsberg manifold, according to \eqref{a-1} and \eqref{a-2}. The second part of the proposition follows analogously.  
\end{proof}

To rigorously define and express a specific Finsler scalar curvature, we first introduce some temporary notations that will be used consistently throughout this manuscript.

The flag curvature tensor is a symmetric 2-tensor. From the Chern Riemannian curvature tensor, we define an asymmetric 2-tensor called the \emph{asymmetric $g$-Ricci curvature}, denoted by $\overline{Ric}:=\bar R_{ij}dx^i\otimes dx^j$ where
\begin{eqnarray}\label{barRij}
	\bar R_{ij}=g^{kl}g_{ih}R_{k\,\,jl}^{\,\,h}.
\end{eqnarray}
Its symmetrization is referred to as the \emph{symmetric $g$-Ricci curvature}, denoted by $\widetilde{Ric}:=\tilde{R}_{ij}dx^i\otimes dx^j$ where
\begin{eqnarray}
	\tilde R_{ij}=\frac{1}{2}(\bar R_{ij}+\bar R_{ji}).
\end{eqnarray}
Furthermore, we define the asymmetric weighted $g$-Ricci curvature is defined by $\overline{Ric}^{\infty}:=\bar R^{\infty}_{ij}dx^i\otimes dx^j$, with 
$$\bar R^{\infty}_{ij}=\bar R_{ij}+\tau_{|i|j},$$ 
and the symmetric weighted $g$-Ricci curvature as $\widetilde{Ric}^{\infty}:=\tilde R^{\infty}_{ij}dx^i\otimes dx^j$, where 
$$\tilde R^{\infty}_{ij}=\tilde R_{ij}+\frac12(\tau_{|i|j}+\tau_{|j|i}).$$
Observe that
\begin{eqnarray}
	Ric(x,y)=\bar R_{ij}y^iy^j=\tilde R_{ij}y^iy^j,\quad\mbox{and}\quad Ric^{\infty}(x,y)=\bar R^{\infty}_{ij}y^iy^j=\tilde R^{\infty}_{ij}y^iy^j.
\end{eqnarray}

We now define the scalar curvature as follows.
\begin{defn}
	For a Finsler manifold $(M,F)$, the scalar curvature is a function on the sphere bundle $SM$ given by
	$$R=\mathrm{tr}_g\widetilde{Ric}=\mathrm{tr}_g\overline{Ric}.$$
	In local coordinates, it can be expressed as
	$$R=g^{ij}R^{\,\,k}_{i\,\,kj}=g^{ij}g^{kl}R_{kijl}.$$
\end{defn}

\section{Lower estimates on gradient Ricci soliton with bounded Ricci curvature}
In this section, we establish lower bounds for the S-curvature and distortion on forward complete Finsler gradient Ricci soliton. Starting from the second variation of arc length, we proceed to prove Theorem \ref{thm-Ricbound}.

\begin{proof}[Proof of Theorem \ref{thm-Ricbound}]
	Let $\gamma_0=\gamma$ be a geodesic on $(M,F)$ parametrized by forward arc length $t$, and let $\gamma_s$ denote a variation of $\gamma$. The second variation formula yields
	\begin{eqnarray}
		L''(0)=-\frac{1}{F(\dot\gamma_0(t))}\int_0^{t_0}g_{\dot\gamma}(\nabla_{\dot\gamma}\nabla_{\dot\gamma}v+R(v),v)dt,
	\end{eqnarray}
	where $R$ denote the flag curvature operator,
	The second variation formula yields $v$ satisfies $g_{\dot\gamma}(\dot\gamma, v)=0$ (see [Ch.5 in \cite{ShenShen2016}]). Here the dot notation $``\,\dot\,\, "$ denotes differentiation along the geodesic with respect to the arc length parameter.
	
	Consider a minimal forward geodesic $\gamma_0$ from $p=\gamma(0)$ to $x=\gamma(t_0)$. Since $L''(0)>0$, we obtain
	\begin{eqnarray}
		\begin{split}
			-\int_0^{t_0}g_{\dot\gamma}(R(v),v)dt&\geq \int_0^{t_0}g_{\dot\gamma}(\nabla_{\dot\gamma}\nabla_{\dot\gamma}v,v)dt\\
			&=\int_0^{t_0}\nabla_{\dot\gamma}(g_{\dot\gamma}(\nabla_{\dot\gamma}v,v))dt-\int_0^{t_0}g_{\dot\gamma}(\nabla_{\dot\gamma}v,\nabla_{\dot\gamma}v)dt.
		\end{split}
	\end{eqnarray}
	For variations with fixed endpoints, this simplifies to
	\begin{eqnarray}\label{gR<gvv}
		\int_0^{t_0}g_{\dot\gamma}(R(v),v)dt\leq\int_0^{t_0}g_{\dot\gamma}(\nabla_{\dot\gamma}v,\nabla_{\dot\gamma}v)dt.
	\end{eqnarray} 
	Let $e_1,\dots,e_n$ be an orthonormal basis at $\gamma(0)$ with respect to $g_{\dot\gamma(0)}$, where $e_n=\dot\gamma(0)$. Parallel transport with respect to $g_{\dot\gamma(t)}$ yields orthonormal vector fields $E_1(t),\dots,E_n(t)=\dot\gamma(t)$  along $\gamma$. Substituting $v=f(t)E_i$ into \eqref{gR<gvv} and summing over $i$ give that
	\begin{eqnarray}\label{fRic<n-1f2}
		\int_0^{t_0}f^2(t)Ric(\dot\gamma(t),\dot\gamma(t))dt\leq (n-1)\int_0^{t_0}|\dot f(t)|^2dt,
	\end{eqnarray}
	where $f(t)$ is a continuous function with $f(0)=f(t_0)=0$, and $Ric(\dot\gamma(t),\dot\gamma(t))$ is the Ricci scalar in the direction $\dot\gamma(t)$.
	Now, define $f(t)$ piecewise as
	\begin{eqnarray}
		f(t)= \begin{cases}t, & \text { if }t\in(0,1); \\ 1, & \text { if }t\in[1,t_0-1]; \\ t_0-1,& \text { if }t\in[t_0-1,t_0].\end{cases}
	\end{eqnarray}
	Substituting into \eqref{fRic<n-1f2}, we derive
	\begin{eqnarray}
		\begin{split}
			\int_0^{t_0}Ric(\dot\gamma(t),\dot\gamma(t))dt&=\int_0^{t_0}f^2(t)Ric(\dot\gamma(t),\dot\gamma(t))dt+\int_0^{t_0}(1-f^2(t))Ric(\dot\gamma(t),\dot\gamma(t))dt\\
			&\leq (n-1)\int_0^{t_0}|\dot f(t)|^2dt+\int_0^{t_0}(1-f^2(t))Ric(\dot\gamma(t),\dot\gamma(t))dt\\
			&\leq 2(n-1)+\frac23\left(\max_{SB_p(1)}|Ric|+\max_{SB_{\gamma(t_0)}(\rho)}|Ric|\right).
		\end{split}
	\end{eqnarray}
	Setting $y=\dot\gamma(t)$ in \eqref{eq-EF}, we observe
	\begin{eqnarray}\label{eq-EF-2}
		\nabla_{\dot\gamma(t)}S=\nabla_{\dot\gamma}\nabla_{\dot\gamma}\tau=\sigma(\gamma,\dot\gamma)-Ric(\dot\gamma,\dot\gamma).
	\end{eqnarray}
	Integrating \eqref{eq-EF-2} along $\gamma$ from $0$ to $t_0$ yields
	\begin{eqnarray}
		\begin{split}
			S(\gamma(t_0),\dot\gamma(t_0))-S(\gamma(0),\dot\gamma(0))&=\dot\tau(\gamma(t_0),\dot\gamma(t_0))-\dot\tau(\gamma(0),\dot\gamma(0))\\
			&\geq\frac{t_0}{2}-\int_0^{t_0}Ric(\dot\gamma,\dot\gamma)dt\\
			&\geq \frac{t_0}{2}-2(n-1)-\frac23\left(\max_{SB_p(1)}|Ric|+\max_{SB_{\gamma(t_0)}(\rho)}|Ric|\right),
		\end{split}
	\end{eqnarray}
	where $\rho$ is the local reversibility at $\gamma(t_0)$.
	Assuming the gradient Ricci solition has bounded Ricci curvature, namely, $|Ric|\leq cF^2$ for some constant $c$, we obtain
	\begin{eqnarray}\label{K0}
		\begin{split}
			S(\gamma(t_0),\dot\gamma(t_0))&\geq \frac{t_0}{2}+S(\gamma(0),\dot\gamma(0))-2(n-1)-\frac{4c}{3}\\
			&=\frac12(t_0-K_0),
		\end{split}
	\end{eqnarray}
	where $K_0$ depends only on the dimension $n$, the Ricci curvature bound $c$, the Finsler metric $F$, and the measure on $S_pM$. Since the S-curvature is 1-homogeneous on $TM_0$, so we see that $S/F$ is a function on $SM$. Thus, it follows from \eqref{K0} that
	\begin{eqnarray}
		S(x,y)\geq\frac12(d(p,x)-K_0),
	\end{eqnarray}
	where $K_0$ only depends on $n$, $c$, $F$ and the measure on $S_pM$.
	
	Integraling \eqref{K0} again yields that
	\begin{eqnarray}
		\tau(\gamma(t_0),\dot\gamma(t_0))\geq \frac{1}{4}(t_0-K_0)^2-\frac14K_0^2-\tau(\gamma(0),\dot\gamma(0)),
	\end{eqnarray}
	and similarly,
	\begin{eqnarray}
		\tau(x,y)\geq \frac{1}{4}(d(p,x)-K_0)^2-K'_0,
	\end{eqnarray}
	where $K_0$ and $K_0'$ are two constants only depending on the dimension $n$, the bounds of Ricci curvature $c$, the reversibility $\rho$ as well as the Finsler metric $F$ and the given measure on $S_pM$.
\end{proof}

The bounded Ricci curvature condition is natural in Finsler geometry. Theorem \ref{thm-Ricbound} thus demonstrates how curvature influences the growth of volume measures via the $\infty$-Einstein metric. For asymmetric essential gradient Ricci soliton, more interesting curvature conditions yield analogous results, which we explore further in Section 7.

\section{A formula on the essential gradient Ricci soliton}
A direct corollary of Theorem \ref{thm-Ricbound} is that the same lower bounds also hold on a forward complete Finsler gradient Ricci soliton, whose $\sigma(x,y)=\frac12$. Furthermore, we can derive both the lower and upper bounds of distortion, S-curvature and the scalar curvature by using a much weaker curvature condition as presented in Theorem \ref{thm-1.2}.
Before showing how the scalar curvature affects the growth of the Riemannian and non-Riemannian geometric quantities in a gradient Ricci soliton, we first derive a key formula crucial for studying the Finsler gradient Ricci soliton.

We denote 
\begin{align}
	\tilde R_{ij}^{\infty}=\frac{1}{2}(\bar R_{ij}+\tau_{|i|j}+\bar R_{ji}+\tau_{|j|i})=\frac{1}{2}(\bar R_{ij}+\bar R_{ji})+\frac{1}{2}(\tau_{|i|j}+\tau_{|j|i})=\widetilde{R}_{ij}+\widetilde{\tau}_{ij}.
\end{align}
Corresponding to that
\begin{align}	
	&\begin{cases}
		\bar{R}_{ij}-\bar{R}_{ji}=-2C_{js}^{t}R_{\,\,ti}^{s}+2C^{t}_{is}R^{s}_{\,\,tj}-I_{s}R^{s}_{\,\,ij}\\
		\tau_{|i|j}-\tau_{|j|i}=I_{s}R_{\,\,ij}^{s},
	\end{cases}
\end{align}
the essential gradient Ricci soliton could be expressed as
\begin{eqnarray}\label{FE}
	\tilde{R}_{ij}^{\infty}=\tilde{R}_{ij}+\tilde{\tau}_{ij}=\bar{R}_{ij}+\tau_{|i|j}+C_{js}^{t}R_{\,\,ti}^{s}-C_{is}^{t}R_{\,\,tj}^{s}=\frac{1}{2}g_{ij}.
\end{eqnarray}
Taking the trace of $i,k$ in the second Bianchi identity 
\begin{align}
	R_{j\,kl|m}^{\,\,i}+R_{j\,lm|k}^{\,\,i}+R_{j\,mk|l}^{\,\,i}=P_{j\,ms}^{\,i}R_{\,\,kl}^{s}+P_{j\,ks}^{\,i}R_{\,\,lm}^{s}+P_{j\,ls}^{\,i}R_{\,\,mk}^{s}
\end{align}
yields
\begin{align}\label{R|3=PR3}
	R_{j\,il|m}^{\,\,i}+R_{j\,lm|i}^{\,\,i}-R_{j\,im|l}^{\,\,i}=P_{j\,ms}^{\,i}R_{\,\,il}^{s}+P_{j\,is}^{\,i}R_{\,\,lm}^{s}+P_{j\,ls}^{\,i}R_{\,\,mi}^{s}.
\end{align}
Contracting \eqref{R|3=PR3} with $g^{jl}$ gives
\begin{align}\label{R|3=PR3-2}
	\bar R^{i}_{\,\,i|m}-\bar R^{i}_{\,\,m|i}-R^{ji}_{\,\,\,\,im|j}=P^{ji}_{\,\,\,\,\,\,ms}R^{s}_{\,\,ij}+P^{ji}_{\,\,\,\,\,is}R^{s}_{\,\,jm}+P^{ji}_{\,\,\,\,\,js}R^{s}_{\,\,mi}.
\end{align}

From the identity
\begin{eqnarray}
	\begin{split}
		R^{ji}_{\,\,\,\,mi|j}=g^{jp}g^{iq}R_{pqmi|j}&=g^{jp}g^{iq}[-R_{qpmi|j}-2(C_{pgs}R^{s}_{\,\,mi})_{|j}] \\
		&=-R^{ij}_{\,\,\,\,mi|j}-2(C^{ji}_{s}R^{s}_{\,\,mi})_{|j},
	\end{split}
\end{eqnarray}
and \eqref{R|3=PR3-2}, we obtain that
\begin{align}
	R_{|m}-2(C^{ij}_sR^{s}_{\,\,mi})_{|j}-P^{ji}_{\,\,\,\,\,ms}R^{s}_{\,\,ij}-(P^{ji}_{\,\,\,\,\,is}-P^{ij}_{\,\,\,\,\,is})R_{\,\,jm}^{s}=2\bar R^{i}_{\,\,m|i}.
\end{align}
Using the relation between the Chern non-Riemannian curvature and the Cartan and Landsberg curvatures,
\begin{eqnarray}
	P_{jims}=-C_{ijs|m}+C_{jms|i}-C_{mis|j}+C_{ij}^{t}L_{tms}-C_{jm}^{t}L_{tis}+C_{mi}^{t}L_{tjs},
\end{eqnarray}
we derive 
\begin{eqnarray}
	\begin{split}
		P^{ji}_{\,\,\,\,\,\,ms}R_{\,\,ij}^{s}&=(C^{j\,\,\,\,|i}_{ms}-C^{i\,\,\,\,\,\,|j}_{ms}-C^{j}_{tm}L^{ti}_s+C^{i}_{tm}L^{tj}_s)R_{ij}^{s}\\
		&=2(C^{j\,\,\,\,|i}_{ms}-C^{j}_{tm}L^{ti}_s)R_{\,\,ij}^{s}.
	\end{split}
\end{eqnarray}
Thus,
\begin{align}\label{2}
	R_{|m}-2(C_{s}^{ij}R_{\,\,mi}^{s})_{|j}-2(C^{j\,\,\,\,|i}_{ms}-C^{j}_{tm}L^{ti}_s)R_{\,\,ij}^{s}-(P_{\,\,\,\,\,is}^{ji}-P_{\,\,\,\,\,is}^{ij})R_{\,\,jm}^{s}=2\bar R^{i}_{\,\,m|i}.
\end{align}
The essential gradient Ricci soliton equation can be written as
\begin{eqnarray}
	\bar{R}^{i}_{\,\,j}+\tau^{|i}_{\,\,\,\,|j}+C_{js}^{t}R_{\,\,t}^{s\,\,i}-C_{s}^{ti}R_{\,\,tj}^{s}=\frac{1}{2}\delta_{j}^{i}.
\end{eqnarray}
Substituting this into \eqref{2} gives 
\begin{eqnarray}\label{RCPR=tCR}
	\begin{split}
		R_{|m}-2(C_{s}^{ij}R^{s}_{\,\,mi})_{|j}&-2(C_{ms}^{j\,\,\,\,|i}-C_{tm}^{j}L_{s}^{ti})R_{\,\,ij}^{s}\\
		&-(P_{\,\,\,\,is}^{ji}-P^{ij}_{\,\,\,\,is})R_{\,\,jm}^{s}=2[-\tau^{|i}_{\,\,\,\,|m|i}-(C_{ms}^{t}R_{\,\,t}^{s\,\,i}-C^{ti}_{s}R_{\,\,tm}^{s})_{|i}].
	\end{split}
\end{eqnarray}
By the Ricci identity \cite{Shen2018}, we find
\begin{eqnarray}
	\begin{split}
		\tau^{|i}_{\,\,\,\,|m|i}&=g^{ik}\tau_{|k|m|i}=g^{ik}(\tau_{|k|i|m}+R_{k\,\,mi}^{\,\,s}\tau_{|s}+R^{s}_{\,\,mi}\tau_{|k;s}) \\
		&=\tau^{|i}_{\,\,\,\,|i|m}+\bar{R}^{s}_{\,\,m}\tau_{|s}+g^{ik}R^{s}_{\,\,mi}\tau_{|k;s} \\
		&=(\frac{n}{2}-R)_{|m}+\bar{R}^{s}_{\,\,m}\tau_{|s}+g^{ik}R^{s}_{\,\,mi}\tau_{k;s}. 
	\end{split}
\end{eqnarray}
Substituting this back into \eqref{RCPR=tCR} yields
\begin{eqnarray}\label{Rt=CRPR}
	\begin{split}
		R_{|m}-2\bar{R}^{s}_{\,\,m}\tau_{|s}=&2g^{ik}R^{s}_{\,\,mi}\tau_{|k;s}+2(C_{ms}^{t}R^{s\,\,i}_{\,\,t}-C_{s}^{ti}R^{s}_{\,\,tm})_{|i}-2(C^{ij}_{s}R^{s}_{\,\,mi})_{|j}\\&-2(C^{j\,\,\,\,\,|i}_{ms}-C^{j}_{tm}L^{ti}_s)R^{s}_{\,\,ij}-(P^{ji}_{\,\,\,\,is}-P^{ij}_{\,\,\,\,is})R^{s}_{\,\,im}.
	\end{split}
\end{eqnarray}
Contracting 
\begin{eqnarray}
	P_{jikl}+P_{ijkl}-2C_{ijs}L^{s}_{kl}+2C_{ijl|k}=0
\end{eqnarray}
with $g^{ij}$ gives
\begin{eqnarray}
	2P_{i\,\,kl}^{\,i}-2I_{s}L_{kl}^{s}+2I_{l|k}=0,
\end{eqnarray}
which, combined with
\begin{eqnarray}
	\tau_{|k;s}=I_{s|k}-I_{t}L_{ks}^{t},
\end{eqnarray}
implies
\begin{eqnarray}
	P_{i\,\,kl}^{\,i}+\tau_{|k;l}=0.
\end{eqnarray}
This leads to
\begin{eqnarray}
	g^{ik}R^{s}_{\,\,mi}\tau_{|k;s}=g^{ik}R^{s}_{\,\,mi}P_{j,\,ks}^{\,j}=g^{jk}R_{\,\,jm}^{s}P_{k\,\,is}^{\,\,i}=P^{ji}_{\,\,\,\,is}R_{\,\,jm}^{s}.
\end{eqnarray}
Replacing the corresponding term in \eqref{Rt=CRPR}, we obtain
\begin{eqnarray}
	\begin{split}
		R_{|m}-2\bar{R}_{\,\,m}^{s}\tau_{|s}=&2(C_{ms}^{t}R^{s\,\,i}_{\,\,t}-C^{ti}_sR_{\,\,tm}^{s})_{|i}-2(C^{ij}_sR^{s}_{\,\,mi})_{|j}\\
		&-2(C^{j\,\,\,\,|i}_{ms}-C^{j}_{tm}L^{ti}_s)R_{\,\,ij}^{s}+(P^{ji}_{\,\,\,\,is}+P^{ij}_{\,\,\,\,is})R_{\,\,jm}^{s}.
	\end{split}
\end{eqnarray}
Noting again that 
\begin{align}
	P^{ji}_{\,\,\,\,is}+P^{ij}_{\,\,\,\,is}=2C^{ij}_t+L^t_{is}-2C^{ij}_{s|i},
\end{align}
we derive
\begin{eqnarray}\label{Rt=CR|}
	\begin{split}
		\frac{1}{2}R_{|m}-\bar R_{\,\,m}^{s}\tau_{|s}=&(C^{t}_{ms}R_{\,\,ti}^{s}-C^{t}_{is}R_{\,\,tm}^{s})^{|i}+(C^{t}_{is}R^{s}_{\,\,tm})^{|i}\\
		&+(C^{t\,\,\,\,|i}_{ms}-C^{t}_{ml}L^{li}_{s})R_{\,\,ti}^{s}+(C^{t}_{il}L^{li}_{s}-C^{t\,\,|i}_{is})R_{\,\,tm}^{s}\\
		=&(C^{t}_{ms}R_{\,\,ti}^{s})^{|i}+(C^{t\,\,\,\,|i}_{ms}-C^{t}_{ml}L^{li}_{s})R_{\,\,ti}^{s}+(C^{t}_{il}L^{li}_{s}-C^{t\,\,|i}_{is})R_{\,\,tm}^{s}.
	\end{split}	
\end{eqnarray}
Contracting \eqref{FE} with $\tau^{|i}$ gives
\begin{eqnarray}
	\tau^{|i}\bar{R}_{im}+\frac{1}{2}(F^2_y(\nabla\tau))_{|m}+\tau^{|i}(C_{ms}^{t}R_{\,\,ti}^{s}-C_{is}^{t}R^{s}_{\,\,tm})=\frac{1}{2}\tau_{|m},
\end{eqnarray}
where $F_y(\nabla\tau)$ denotes the norm of the horizontal Chern-covariant derivative of the distorsion with respect to $y$. Combining this with \eqref{Rt=CR|} yields the following key formula.

\begin{lem}
	Let $(M,F)$ be an essential gradient Ricci soliton satisfying \eqref{def-VV}, with distortion $\tau$ and scalar curvature $R$. Then
	\begin{eqnarray}\label{Rtau=C}
		\frac12(R+F^2_y(\nabla\tau)-\tau)_{|m}+\tau^{|i}(C_{ms}^{t}R_{\,\,ti}^{s}-C_{is}^{t}R^{s}_{\,\,tm})=K_m,
	\end{eqnarray}
	where $K_m=(C^{t}_{ms}R_{\,\,ti}^{s})^{|i}+(C^{t\,\,\,\,|i}_{ms}-C^{t}_{ml}L^{li}_{s})R_{\,\,ti}^{s}+(C^{t}_{il}L^{li}_{s}-C^{t\,\,|i}_{is})R_{\,\,tm}^{s}$.
\end{lem}

From \eqref{Rtau=C}, we also deduce that
\begin{eqnarray}\label{Rtau=C0}
	\frac12(R+F^2_y(\nabla\tau)-\tau)_{|0}-\tau^{|i}C_{is}^{t}R^{s}_{\,\,t}=(C^{t}_{il}L^{li}_{s}-C^{t\,\,|i}_{is})R_{\,\,t}^{s}.
\end{eqnarray}

\section{Some special cases}
In this section, we examine the fundamental equation \eqref{Rtau=C} on specific classes of Finsler manifolds.
\subsection{On Berwald manifolds}	
For a Berwald metric, \eqref{Rtau=C} simplifies to
\begin{eqnarray}\label{Rtau=CB}
	\frac12(R+F^2_y(\nabla\tau)-\tau)_{|m}+\tau^{|i}(C_{ms}^{t}R_{\,\,ti}^{s}-C_{is}^{t}R^{s}_{\,\,tm})=C^{t}_{ms}R_{\,\,ti}^{s\,\,\,|i}.
\end{eqnarray}
Similarly, \eqref{Rtau=C0} reduces to
\begin{eqnarray}
	\frac12(R+F^2_y(\nabla\tau)-\tau)_{|0}=\tau^{|i}C_{is}^{t}R^{s}_{\,\,t}.
\end{eqnarray}

Since a Berwald metric with Busemann-Hausdorff volume form necessarily has vanishing S-curvature and constant volume function $\sigma_B$, equations \eqref{Rtau=CB} and \eqref{Rtau=C0} further simplify to
\begin{eqnarray}\label{R|0=0}
	R_{|m}=2C^{t}_{ms}R_{\,\,ti}^{s\,\,\,|i}, \quad\mbox{and}\quad R_{|0}=0.
\end{eqnarray}
This implies that the scalar curvature $R$ remains constant along any geodesic on the manifold. Assuming the gradient Ricci soliton is forward complete, $R$ becomes a 0-homogeneous function dependent solely on $y$. 
Consequently, we derive the following result.
\begin{thm}
	Let $(M,F,\mu)$ be a Berwald essential Finsler gradient Ricci soliton equipped with the Busemann-Hausdorff volume form. The scalar curvature $R$ is a function defined on $S_{x_0}M$, where $x_0$ is an arbitrary point on $M$ and $S_{x_0}M$ denotes the indicatrix at $x_0$. As a result, $R$ is a bounded function.
\end{thm}

\subsection{On Finsler manifolds with isotropic S-curvature}
Einstein Finsler manifolds with isotropic S-curvature have been extensively studied in \cite{ShenSun2021,ShenZhao2022,LMW2024,Xia2024}, among others. For a metric admitting isotropic S-curvature, i.e., $S=(n+1)c(x)F$, where $c(x)$ is a scalar function on $M$, the distortion can be expressed as $\tau_{|i}=(n+1)c(x)l_i$, with $l_i=g_{ij}\frac{y^j}{F}$, since this choice preserves the S-curvature of the original manifold. This yields
\begin{align}
	F^2_y(\nabla\tau)=(n+1)^2c(x),
\end{align}
and 
\begin{eqnarray}\label{Rc=CR}
	\frac12(R_{|0}+(n+1)^2(c^2)_{|0}-(n+1)cF)=(C^{t}_{il}L^{li}_{s}-C^{t\,\,|i}_{is})R_{\,\,t}^{s}.
\end{eqnarray}
Assuming $R_{|0}\geq(n+1)K_2F$ and $(C^{t}_{il}L^{li}_{s}-C^{t\,\,|i}_{is})R_{\,\,t}^{s}\leq \frac{(n+1)}{2}K_1F$ for positive constants $K_1$ and $K_2$, we obtain 
\begin{eqnarray}
	(c^2)_{|0}\leq \frac{1}{n+1}(K_1-K_2+c)F,
\end{eqnarray}
which is equivalent to
\begin{eqnarray}
	\frac{2(n+1)cc_{|0}}{K_1-K_2+c}\leq F.
\end{eqnarray}
One can deduce from it that
\begin{eqnarray}\label{c|0<F}
	[c-(K_1-K_2)\log(K_1-K_2+c)]_{|0}\leq \frac{F}{2(n+1)}.
\end{eqnarray}
Hence for any two points $x_1$ and $x_2$ on $M$, there is a minimal geodesic $L$ from $x_1$ to $x_2$, whose length gives the forward distance $d(x_1,x_2)$ from $x_1$ to $x_2$. Integrating \eqref{c|0<F} along $L$ gives that
\begin{eqnarray}
	c(x_2)-c(x_1)-(K_1-K_2)\log\frac{K_1-K_2+c(x_2)}{K_1-K_2+c(x_1)}\leq\frac{1}{2(n+1)}d(x_1,x_2).
\end{eqnarray}
This leads to the following estimate for the S-curvature.
\begin{thm}
	Let $(M,F,p)$ be a forward complete essential Finsler gradient Ricci soliton with a pole $p$, admitting isotropic S-curvature $S=(n+1)c(x)F$ for a scalar function $c(x)$. Suppose the curvature satisfies $R_{|0}\geq-(n+1)K_2F$ and $(C^{t}_{il}L^{li}_{s}-C^{t\,\,|i}_{is})R_{\,\,t}^{s}\leq \frac{(n+1)}{2}K_1F$ for positive constants $K_1$ and $K_2$. Then for any point $x$ on $M$, 
	\begin{eqnarray}
		c(x)-(K_1+K_2)\log(K_1+K_2+c(x))\leq\frac{1}{2(n+1)}d(p,x)+K_3,
	\end{eqnarray}
	where $K_3=c(p)-(K_1+K_2)\log(K_1+K_2+c(p))$ is a constant, and $d(p,x)$ denotes the forward distance from $p$ to $x$. 
\end{thm}

\begin{rem}
	\begin{itemize}
		\item The curvature condition $R_{|0}\geq-(n+1)K_2F$ is a mild constraint on the scalar curvature $R$. essentially requiring that $R$ exhibits a linear growth lower bound of the form $R\geq C-(n+1)K_2d(p,x)$.
		\item The curvature condition $(C^{t}_{il}L^{li}_{s}-C^{t\,\,|i}_{is})R_{\,\,t}^{s}\leq \frac{(n+1)}{2}K_1F$ is reasonable, since it is automatically satisfied on compact manifolds or precompact sets. More special examples could be referred to \cite{Xia2024}.
	\end{itemize}
	
\end{rem}

A special case is that the metric admits a constant S-curvature, namely, $c(x)$ is a constant $c$. Then \eqref{Rc=CR} implies that
\begin{eqnarray}\label{R0<KcF}
	R_{|0}\leq (n+1)(K_1+c)F,
\end{eqnarray}
provided $(C^{t}_{il}L^{li}_{s}-C^{t\,\,|i}_{is})R_{\,\,t}^{s}\leq \frac{(n+1)}{2}K_1F$. Integrating along any forward geodesic originating from $p$ yields the following corollary.
\begin{cor}
	Let $(M,F,p)$ be a forward complete essential Finsler gradient Ricci soliton with a pole $p$and constant S-curvature $S=(n+1)cF$. If the curvature bound $(C^{t}_{il}L^{li}_{s}-C^{t\,\,|i}_{is})R_{\,\,t}^{s}\leq \frac{(n+1)}{2}K_1F$ holds for a positive constant $K_1$, then for any $x$ on $M$, it satisfies that
	\begin{eqnarray}
		R\leq(n+1)(K_1+c)d(p,x)+K_4,
	\end{eqnarray}
	where $K_4$ is a constant determined by $S_pM$. 
\end{cor}
\begin{proof}
	For any point $x$ on $M$, there is a minimal forward geodesic $\sigma$ from $p$ to $x$, for the forward completeness of $M$. Integrating \eqref{R0<KcF} along $\sigma$ yields
	\begin{eqnarray}\label{R-R<Kcd}
		R(x,y_x)-R(p,y_p)\leq (n+1)(K_1+c)d(p,x),
	\end{eqnarray}
	where $y_x\in T_xM$ and $y_p\in T_pM$ is the unique vector parallel-transported to $y_x$ along $\sigma$. Thus, for any $(x,y)\in SM$, it follows from \eqref{R-R<Kcd} that 
	\begin{eqnarray}
		R(x,y)\leq (n+1)(K_1+c)d(p,x)+ K_4,
	\end{eqnarray}
	with $K_4=\sup_{y\in S_pM}R(p,y)$.
\end{proof}

\section{On asymmetric and strongly asymmetric essential gradient Ricci solitons}
The symmetry of the Ricci-type curvature $\bar R$ renders the Einstein equation more natural, though at the cost of losing some tangent space information. In this section, we focus on two sub-concepts to derive an upper bound estimate for the scalar curvature using the forward distance function and its lower bound.
According to Definition \ref{def-asymFE}, a Finsler metric is called an asymmetric essential gradient Ricci soliton, if it satisfies that
\begin{eqnarray}
	\bar{R}_{ij}+\tau_{|i|j}=\frac{1}{2}g_{ij}.
\end{eqnarray}
This condition is equivalent to requiring $C_{ms}^{t}R_{\,\,ti}^{s}-C_{is}^{t}R^{s}_{\,\,tm}=0$. Combining these, we obtain
\begin{eqnarray}\label{Rtau=Km}
	\frac12(R+F^2_y(\nabla\tau)-\tau)_{|m}=K_m,
\end{eqnarray}
where $K_m=(C^{t}_{ms}R_{\,\,ti}^{s})^{|i}+(C^{t\,\,\,\,|i}_{ms}-C^{t}_{ml}L^{li}_{s})R_{\,\,ti}^{s}+(C^{t}_{il}L^{li}_{s}-C^{t\,\,|i}_{is})R_{\,\,tm}^{s}$.

We now analyze \eqref{Rtau=Km} to estimate the scalar curvature, distortion, and S-curvature. The method applies uniformly to both asymmetric and strongly asymmetric cases, differing only in the expression of
$K_m$ in \eqref{Rtau=Km}. 

\begin{thm}\label{thm-tauSR}
	Let $(M,F,p)$ be a forward complete asymmetric essential Finsler gradient Ricci soliton with a pole $p$, where the scalar curvature satisfies the linear growth condition $R\geq \gamma d(p,x)-\alpha$ for positive constants $\alpha$ and $\gamma$. If the curvature bound $(C^{t}_{il}L^{li}_{s}-C^{t\,\,|i}_{is})R_{\,\,t}^{s}\leq \frac{(n+1)}{2}K_1F$ holds with $K_1\leq (n+1)\gamma$, then the following upper bounds hold for distortion, S-curvature, and scalar curvature.
	\begin{align*}
		\tau(x,y)&\leq\frac14\left[d(p,x)+K_5\right]^2-\alpha-\beta,\\
		|S|&\leq\left[\frac12 d(p,x)+K_5\right]F,\\
		R&\leq \frac14\left[d(p,x)+K_5\right]^2+\gamma d(p,x)-\alpha,
	\end{align*}
	where $\beta$ and $K_5=2\sqrt{\sup_{y\in S_pM}\tau(p,y)+\alpha+\beta}$ are two constants depending only on $S_pM$.
\end{thm}

\begin{proof}
	Rewriting \eqref{Rtau=Km} globally as
	\begin{eqnarray}\label{Rtau=K0}
		\frac1{2}(R+F^2_y(\nabla\tau)-\tau)_{|0}=K_0,
	\end{eqnarray}
	where $K_0=(C^{t}_{il}L^{li}_{s}-C^{t\,\,|i}_{is})R_{\,\,t}^{s}$. Under the given curvature conditions, $K_0\leq \frac12\gamma F$, which is locally satisfied since $K_0$ is a 1-homogeneous function on the tangent bundle. Therefore, we have
	\begin{eqnarray}\label{Rtau<F}
		(R+F^2_y(\nabla\tau)-\tau)_{|0}\leq \gamma F.
	\end{eqnarray}
	SBy forward completeness, for any  $x\in M$, there exists a forward geodesic $\sigma$ from the given fixed point $p$ to $x$. Integrating \eqref{Rtau<F} along $\sigma$ yields
	\begin{eqnarray}\label{Rtau<gammadf}
		R+F^2_y(\nabla\tau)-\tau\leq \gamma d(p,x)+f(y),
	\end{eqnarray}
	where $f(y)$ is a 0-homogeneous function on $S_pM$, uniformly bounded as $f(y)\leq |\sup_{S_{p}M}f(y)|=:\beta$. Combining this with the linear growth assumption on $R$, we can deduce that
	\begin{eqnarray}\label{nabtau-tau<beta}
		F^2_y(\nabla\tau)-\tau\leq\alpha+\beta.
	\end{eqnarray}
	Let $\chi=\tau+\alpha+\beta$. Then \eqref{nabtau-tau<beta} is equal to 
	\begin{eqnarray}
		|\nabla\sqrt{\chi}|\leq \frac12.
	\end{eqnarray}
	Since $\sqrt{\chi}$ is Lipschitz on $SM$, we could deduce from it that
	\begin{eqnarray}\label{sqrtchi<d}
		|\sqrt{\chi(x,y_x)}-\sqrt{\chi(p,y_p)}|\leq \frac12d(p,x),
	\end{eqnarray}
	for any $y_x\in T_xM$ and $y_p$ is the unique vector on $T_pM$ such that $y_x$ is the parallel transport of $y_p$ along $\sigma$. So it could be deduced that
	\begin{eqnarray}
		\tau(x,y)\leq\frac14\left[d(p,x)+2\sqrt{\sup_{y\in S_pM}\tau(p,y)+\alpha+\beta}\right]^2-\alpha-\beta.
	\end{eqnarray}
	From \eqref{nabtau-tau<beta}, we also derive
	\begin{eqnarray}
		F_y(\nabla\tau)\leq\frac12 d(p,x)+\sqrt{\sup_{y\in S_pM}\tau(p,y)+\alpha+\beta},
	\end{eqnarray}
	so that 
	\begin{eqnarray}
		|S|\leq \left[\frac12 d(p,x)+\sqrt{\sup_{y\in S_pM}\tau(p,y)+\alpha+\beta}\right]F.
	\end{eqnarray}
	Moreover, it follows from \eqref{Rtau<gammadf} that
	\begin{eqnarray}
		R\leq \frac14\left[d(p,x)+2\sqrt{\sup_{y\in S_pM}\tau(p,y)+\alpha+\beta}\right]^2+\gamma d(p,x)-\alpha.
	\end{eqnarray}
\end{proof}

This method originates from the work of Perelman \cite{Pere2003} and was refined by H. Cao and D. Zhou \cite{CZ2010}. 

Not surprisingly, the $\infty$-Einstein Finsler metric connected the pure curvatures and the distortion, allowing the scalar curvature bound to depend on distortion. If we choose the volume function to endure the uniform finiteness of the distortion on the whole asymmetric essential gradient Ricci soliton, that is, $\tau\leq \delta$ on $SM$ for some constant $\delta$, \eqref{Rtau<gammadf} provides a linear upper bound that
\begin{eqnarray}
	R\leq \gamma d(p,x) + \beta+ \delta,
\end{eqnarray}
where $\beta=\sup_{S_{x_0}M}f$ is a constant depending only on $SM$. Theorem \ref{thm-tauSR} shows how the scalar curvature's lower bound influences its own upper bound and the distortion on a Finsler gradient Ricci soliton. 


One may not be satisfied with the linear growth condition of the scalar curvature $R$ and may think that a constant lower bound will be better. In order to do so, more non-Riemannian curvature conditions need to be constrained. One of this is to suppose that $K_m$ is a horizontal derivative of a bounded function $h(x,y)$ on the sphere bundle. That is, 
\begin{eqnarray}\label{Km=hm}
	K_m=\frac12h_{|m},
\end{eqnarray}
with $h(x,y)\leq \gamma$ being bounded from above by a positive constant.
Plugging \eqref{Km=hm} into \eqref{Rtau=Km} yields that
\begin{eqnarray}
	R+F^2_y(\nabla\tau)-\tau=h(x,y)+f(y),
\end{eqnarray}
for some function $f(y)$ defined only on the tangent sphere $S_{x_0}M$ at any point $x_0$, with $f(y)\leq |\sup_{S_{x_0}M}f|=:\beta$. 
By the same argument in the proof of Theorem \ref{thm-tauSR}, we immediately obtain the the following corollary.

\begin{cor}
	Let $(M,F,p)$ be a forward complete asymmetric essential Finsler gradient Ricci soliton with a pole $p$, whose scalar curvature is bounded from below by $R\geq-\alpha$ for some constant $\alpha>0$. If $(C^{t}_{il}L^{li}_{s}-C^{t\,\,|i}_{is})R_{\,\,t}^{s}=h_{|0}$, with $h=h(x,y)\leq \gamma$ a bounded 0-homogeneous function on $SM$. Then we have 
	\begin{align*}
		\tau(x,y)&\leq\frac14\left[d(p,x)+K_5\right]^2-\alpha-\beta-\gamma\\
		|S|&\leq\left[\frac12 d(p,x)+K_5\right]F,\\
		R&\leq \frac14\left[d(p,x)+K_5\right]^2-\alpha,
	\end{align*}
	where $\beta$ and $K_5=2\sqrt{\sup_{y\in S_pM}\tau(p,y)+\alpha+\beta+\gamma}$ are two constants depending only on $S_pM$.
\end{cor} 

\begin{rem}
	The condition of \eqref{Km=hm} is restrictive, as $h$ must satisfy a system of differential equations on $SM$, which holds trivially for flat manifolds but is challenging in general.
\end{rem}

Inspired by the integration by parts method as did in \cite{FMZ2008}, we can improve the result in Section 4 by using weaker curvature conditions. Specifically, we get the following result.
\begin{thm}\label{thm-lowbound}
	Under the same hypotheses as Theorem \ref{thm-tauSR}, with $|(C^{t}_{il}L^{li}_{s}-C^{t\,\,|i}_{is})R_{\,\,t}^{s}|\leq \frac{(n+1)}{2}K_1F$, the following lower bounds hold.
	\begin{align*}
		\tau(x,y)&\geq\frac14\left[d(p,x)-K_6\right]^2-\alpha-\beta,\\
		R&\geq \frac14\left[d(p,x)+K_7\right]^2+\gamma d(p,x)+\delta-\alpha-\beta,
	\end{align*}
	where $\beta$, $\delta$ depand on $S_pM$, and $K_6$, $K_7$ depand on $n$, $\gamma$, $F$ and the volume form on the local tangent sphere bundle $SB_p(1)$.
\end{thm}

\begin{proof}
	Integrating \eqref{eq-EF-2} again along the geodesic $\sigma$ from $t=1$ to $t=t_0-1$ that
	\begin{eqnarray}\label{ss>cRicint}
		\begin{split}
			S(\sigma(t_0-1)&,\dot\sigma(t_0-1))-S(\sigma(1),\dot\sigma(1))=\int_1^{t_0-1}\dot\tau(\sigma(t),\dot\sigma(t))dt\\
			=&\int_{1}^{t_0-1}\left[\frac12-Ric(\sigma(t),\dot\sigma(t))\right]dt\\
			=&\frac{t_0-2}2-\int_{1}^{t_0-1}f^2(t)Ric(\sigma(t),\dot\sigma(t))dt\\
			\geq&\frac{t_0-2}2-2(n-1)-\frac13\max_{SB_{p}(1)}|Ric|+\int_{t_0-1}^{t_0}f^2(t)Ric(\sigma(t),\dot\sigma(t))dt.
		\end{split}
	\end{eqnarray}
	Adopting the intergration by parts method, one may find the equivalent form of the last term on the RHS of \eqref{ss>cRicint} as
	\begin{eqnarray}\label{fRic=ctauint}
		\begin{split}
			\int_{t_0-1}^{t_0}f^2(t)Ric(\sigma(t),\dot\sigma(t))dt=&\int_{t_0-1}^{t_0}f^2(t)(\frac12-\nabla_{\dot\sigma}S)dt\\
			=&\frac12\int_{t_0-1}^{t_0}f^2(t)dt-\int_{t_0-1}^{t_0}f^2(t)\nabla_{\dot\sigma}S dt\\
			=&\frac{1}6+S(\sigma(t_0-1),\dot\sigma(t_0-1))-2\int_{t_0-1}^{t_0}f(t)S dt.
		\end{split}
	\end{eqnarray}
	Plugging \eqref{fRic=ctauint} into \eqref{ss>cRicint} makes that
	\begin{eqnarray}\label{ftau>c-Ric+S}
		\begin{split}
			2\int_{t_0-1}^{t_0}f(t)S dt\geq\frac{t_0-2}{2}-2(n-1)-\frac13\max_{SB_p(1)}|Ric|+\frac16+S(\sigma(1),\dot\sigma(1)).
		\end{split}
	\end{eqnarray}
	According to the curvature conditions, we deduce from Theorem \ref{thm-tauSR} that
	\begin{eqnarray}
		\begin{split}
			\max_{t_0-1\leq t\leq t_0}|S(\sigma(t),\dot\sigma(t))|&\leq \max_{t_0-1\leq t\leq t_0}\sqrt{\tau(\sigma(t),\dot\sigma(t))+\alpha+\beta}\\
			&\leq \sqrt{\tau(\sigma(t_0),\dot\sigma(t_0))+\alpha+\beta}+\frac{1}{2},
		\end{split}		
	\end{eqnarray}
	where we have utilized \eqref{sqrtchi<d} at the last inequality. Thus, the LHS of \eqref{ftau>c-Ric+S} could be estimated by
	\begin{eqnarray}
		\begin{split}
			2\int_{t_0-1}^{t_0}f(t)S dt&\leq \max_{t_0-1\leq t\leq t_0}|S(\sigma(t),\dot\sigma(t))|\cdot 2\int_{t_0-1}^{t_0}f(t) dt\\
			&\leq \sqrt{\tau(\sigma(t_0),\dot\sigma(t_0))+\alpha+\beta}+\frac{1}{2},
		\end{split}
	\end{eqnarray}
	and hence,
	\begin{eqnarray}
		\begin{split}
			\sqrt{\tau(\sigma(t_0),\dot\sigma(t_0))+\alpha+\beta}&\geq \frac{t_0}{2}-2n+\frac76-\frac13\max_{SB_p(1)}|Ric|+S(\sigma(1),\dot\sigma(1))\\
			&\geq \frac12(t_0-K_6),
		\end{split}
	\end{eqnarray}
	where $K_6$ depends only on the dimension $n$ and the Finsler metric $F$ as well as the volume form on the local tangent sphere bundle $SB_p(1)$. Therefore, by the same arguement as in the proof of Theorem \ref{thm-tauSR}, we can get that
	\begin{eqnarray}
		\sqrt{\tau(x,y)+\alpha+\beta}\geq \frac{1}{2}(d(p,x)-K_6),
	\end{eqnarray}
	where $K_6$ is a constant depending only on $n$, $F$ and the volume form on the local tangent sphere bundle $SB_p(1)$.
	
	It follows from \eqref{Rtau=K0} and the curvature conditions that
	\begin{eqnarray}
		R+F_y^2(\nabla \tau)-\tau\geq -\gamma d(p,x)+f(y),
	\end{eqnarray}
	with $f(y)\geq\inf_{s_pM}f(y)=:\delta$. Hence it satisfies that
	\begin{eqnarray}
		\begin{split}
			R&\geq \tau-\gamma d(p,x)+\delta\\
			&\geq \frac14\left[d(p,x)-K_6\right]^2-\gamma d(p,x)+\delta-\alpha-\beta\\
			&=\frac14\left[d(p,x)-K_7\right]^2+K_8,
		\end{split}		
	\end{eqnarray}	
	where $K_7$ is a constant depending on $n$, $\gamma$, $F$ as well as the volume form on the local tangent sphere bundle $SB_p(1)$, and $K_8$ is a constant depending on $\alpha$, $\beta$ and $\delta$.
\end{proof}

Theorem \ref{thm-1.2} is a combination of Theorems \ref{thm-tauSR} and \ref{thm-lowbound} by unifying constants $\delta=\beta$.

    
    



	\hskip -0.6cm
	Bin Shen\\
	School of Mathematics, Southeast University, Nanjing 211189, P. R. China\\
	E-mail: shenbin@seu.edu.cn

\end{document}